\documentclass[10pt]{article}
\usepackage{amssymb}
\usepackage{amsfonts}
\usepackage{amsmath}
\usepackage{graphicx}
\usepackage{theorem}
\usepackage[figurename=--\ Figure,labelsep=endash,labelfont=bf]{caption}
\usepackage{hyperref}
\hypersetup{
    colorlinks=true,
    linkcolor=blue,
    filecolor=magenta,      
    urlcolor=cyan,
    pdftitle={Overleaf Example},
    }

\newtheorem{theorem}{Theorem}

\newtheorem{corollary}{Corollary}

\begin{document}
\title{On the non-existence of certain real algebraic surfaces}
\author{Miguel Angel Guadarrama-Garc\'{i}a }
\date{\today}

\maketitle

\begin{abstract}
In this note is given an algebraic solution to the problem 1997-6 proposed by D. A. Panov in the list of Arnold's problems \cite{Arnld2b}. In particular, it is shown that there does not exist a real polynomial function $f$ on the real euclidean plane, whose Hessian is positive in an open set bordered by smooth connected curve, and the parabolic curve of the graph of $f$ has only one special parabolic point with index $+1$. Besides, we find conditions on $f$ so that its graph has more special parabolic points with index -1 than with index +1.
\end{abstract}

\section{The Panov Problem and its Solutions}

Problem 1997-6 in Arnold's problems \cite{Arnld2b} is attributed to D. A. Panov. This problem deals with the affine geometric structure of generic surfaces immersed in the real three-dimensional affine space, which appears by the contact of surfaces with straight lines.  It asks the following.

\emph{Does a generic function $f$ exist on the Euclidean plane, whose Hessian is positive in a region, bordered by smooth connected curve, and the field of asymptotic directions on the parabolic curve of the graph of $f$ has only one special parabolic point with index $+1$? Is it true that the number of special parabolic points with index $-1$ on such curve is not less than the number of special parabolic points with index $+1$?}

The special parabolic points of a smooth surface $S\subset \mathbb{R}^{3}$ are the points of tangency of the asymptotic direction  with the parabolic curve of $S$ (see section \ref{genericsrfcs} for an alternative definition). Over the subset $U$ of hyperbolic points of $S$ (where de Gaussian curvature of $S$ is negative), the bivalued field of asymptotic directions defines a two-sheeted covering surface $M$ in the manifold of the non-oriented tangent elements  $U \times \mathbb{R}P^{1}$. Each hyperbolic point is lifted to the two asymptotic directions at that point.

For generic surfaces, the surface $M$ is smoothly continued by the asymptotic directions at the parabolic points of $S$. The critical line $\Sigma \subset M$ of the canonical projection $\pi : M \longrightarrow S$ 
lies above the parabolic curve.

The asymptotic directions at the hyperbolic points of $S$ are lifted to a direction field on the surface $M$. This field of directions on the surface $M$ is smoothly continued to the critical line $\Sigma $, except for those special parabolic points, where the unique asymptotic direction is tangent to the parabolic curve of $S$.

\begin{figure}[htb]    
\begin{center}
\includegraphics[width=3.77in]{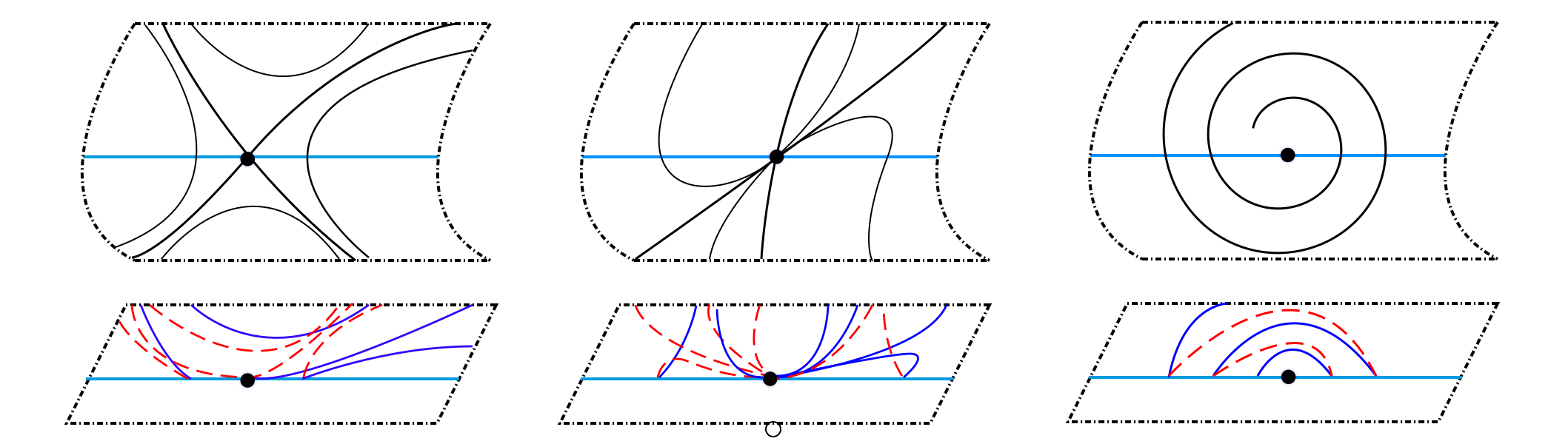}
\caption{Special parabolic points of index -1 (at left) and +1 (at center and right).}
 \label{dblplsindex}
\end{center}
\end{figure}

In general, each special parabolic point is lifted at a singular point of a smooth vector field  $\xi$, which is defined in a neighbourhood of the point in question on the surface $M$. A special parabolic point $p$ can be lifted in a saddle  point (in this case the \emph{index of $p$ is $-1$}),
either a node or focus (in this case the \emph{index of $p$ is $+1$}), 
 topological models of asymptotic curves in a neighbourhood of a special parabolic point and its lifts are shown in Figure \ref{dblplsindex}.

In the case of polynomial functions, we have a negative answer to the first question in the Panov's problem.  Let $$S^{n}_{f}=\left\lbrace (x,y,z)\in \mathbb{R}^{3} \vert z=f(x,y)\right\rbrace $$ be the graph of a polynomial function $f\in \mathbb{R}[x,y]$ of degree \(n\). Write  $$f = \sum _{i=0} ^{n} f_{i},$$ where $f_{i}\in \mathbb{R}[x,y]$ is a homogeneous polynomial of degree $i$.

\begin{corollary}\label{maincor1}
There are no polynomials $f\in \mathbb{R}[x,y]$ of degree \(n\geq 3\), such that all the real factors of $f_{n}$ are simple, $S^{n}_{f}$ is a generic surface, the parabolic curve of  $S^{n}_{f}$ is connected, compact, has one single special parabolic point with index $+1$, and the set of hyperbolic points of $S^{n}_{f}$ is unbounded.
\end{corollary}

It was already known that there are no a generic functions $f$ whose Hessian is negative in a bounded region, bordered by smooth connected closed curve, and such that the field of asymptotic directions on the parabolic curve of the graph of $f$ has only one special parabolic point with index $-1$, see \cite{Blckr&Wlsn}, \cite{Bnchff&Thm} and \cite{Urbcr}. This is because the two-sheeted covering surface $M $ is diffeomorphic to the two dimensional sphere $\mathbb{S}^{2}$. Therefore, if there were a special parabolic point with index $-1$, the vector field  $\xi$ must have others singular points of positive index.

The Corollary \ref{maincor1} is obtained as a special case of the Theorem \ref{mainthrm1}. We identify the parabolic curve of the surface $S^{n}_{f}$ with an affine real algebraic curve called Hessian curve of $f$ (see section \ref{hssncurve}). The closure in $\mathbb{R}P^{2}$ of this curve is the zero locus of a homogeneous polynomial $H_{f} \in \mathbb{R}[x,y,z]$ and is called projective Hessian curve of $f$ (see Section \ref{prhssncurve}). Let $H_{\leq 0}$ be the set of points $[x:y:z]$ in $\mathbb{R}P^{2}$ such that  $H_{f}(x,y,z)\leq 0 $. 

\begin{theorem}\label{mainthrm1}
There are no polynomial functions $f\in \mathbb{R}[x,y]$ of degree \(n\geq 3\), such that all the real factors of $f_{n}$ are simple, $S^{n}_{f}$ is a generic surface with a single special parabolic point of index $+1$, the Hessian projective curve is connected, transverse to the line at infinity, and the set $H_{\leq 0}\subset \mathbb{R}P^{2} $ is not orientable.
\end{theorem}

The case $n=3$ in Theorem \ref{mainthrm1} was previously known: a generic surface $S^{3}_{f}$, with compact parabolic curve, has three special parabolic points with index $-1$, and the set of hyperbolic points of $S^{n}_{f}$ is unbounded. If the parabolic curve is non-compact,  $S^{3}_{f}$ has a special parabolic point with index $-1$ (see \cite{Pnv1} and \cite{H&O&S1}). For $n=1,2$ the set of parabolic points of $S^{n}_{f}$ is the entire surface $S^{n}_{f}$, or empty. In both cases, there are no special parabolic points.  

In contrast, we prove that the second question of the problem of Panov admits an affirmative answer. Let $ \chi (B^{\epsilon} ) $ denote the Euler-Poincar\'{e} characteristic of a smooth bordered surface $B^{\epsilon}, \epsilon\in \left\lbrace -,+\right\rbrace $  (they are specified in section \ref{prhssncurve}). 

\begin{theorem}\label{mainthrm2}
Let $f\in \mathbb{R}[x,y]$ be a polynomial function of degree \(n\geq 3\), such that $S^{n}_{f}$ is a generic surface, and the Hessian projective curve is transverse to the line at infinity. Assume that every real factor of $f_{n}$ is simple and let $k$ be the number of real linear factors of this homogeneous polynomial. If the set $B^{-}$ ( $B^{+}$) contains the set of hyperbolic points of $S^{n}_{f}$, and $k-2\chi (B^{-})\geq 0$ (respectively $k-2\chi (B^{+})\geq 0$), then the number of special parabolic points with index $-1$ is equal to or greater than the number of special parabolic points with index $+1$.
\end{theorem}

Theorems \ref{mainthrm1} and \ref{mainthrm2} are proved using a Poincar\'{e}-Hopf type theorem for surfaces with boundary, which also applies to non-orientable surfaces in the real projective plane (see Theorem \ref{thglobalbminus}). This theorem was obtained by extending the bivalued field of asymptotic directions to  the real projective plane; see section \ref{sec:7} for details. A more complete outline of these ideas can be found in \cite{Gdrrm&Ortz1}. Corollary \ref{maincor1}, theorem \ref{mainthrm1}, and a less general version of theorem \ref{mainthrm2} are contained in the author's doctoral dissertation \cite{Gdrrmtss}.

The ensuing sections offer a brief review of the context surrounding Theorems \ref{mainthrm1} and \ref{mainthrm2}. The proof of these theorems is reserved for the final section of this note.

\section{A classification of the points of a surface $S^{n}_{f}$}\label{genericsrfcs}

Let $f\in \mathbb{R}[x,y]$ be a  polynomial function of degree $n\geq 3$ and $S^{n}_{f}$ be the graph of $f$. A line $l \subset \mathbb{R}^{3}$ is \emph{tangent} at \(p\in S^{n}_{f}\) if the order of contact of $l$ with \(S^{n}_{f}\) at $p$ is equal to or greater than two. In this case, we write  $l\in T_{p}S^{n}_{f}$. We say that a tangent line $l\in T_{p}S^{n}_{f}$ is \emph{simple} if the order of contact at $p$ is two. An \emph{asymptotic line} is a tangent line $l\in T_{p}S^{n}_{f}$ with order of contact at $p$ greater than two. 

A point \(p\in S^{n}_{f}\) is:
\begin{enumerate}
\item  \emph{elliptic} if every tangent line $l\in T_{p}S^{n}_{f}$ is simple.

\item \emph{generic parabolic} if there exists a unique asymptotic line having order of contact at $p$ equals to three.

\item \emph{special parabolic} if there exists a unique asymptotic line having order of contact at $p$ equal to four.

\item \emph{generic hyperbolic} if there are two different asymptotic lines with order of contact at $p$ equal to three.

\item \emph{special hyperbolic} if there are two different asymptotic lines with order of contact at $p$ equal to four.

\item \emph{generic inflection} if there are two different asymptotic lines with order of contact at $p$ equal to three and four, respectively.

\item \emph{bi-inflection} if there are two different asymptotic lines with order of contact at $p$ equal to three and five, respectively.
\end{enumerate}

Fix a natural number $n\geq 3$. Let $\mathbf{S^{n}}$ be the space of all algebraic surfaces $S^{n}_{f}$, with $f\in \mathbb{R}[x,y]$ of degree $n$. The conditions imposed by contact order on each class of points define a positive codimensional submanifold in the space of surfaces $\mathbf{S^{n}}$.  We say that a surface $S^{n}_{f}$ is  \emph{generic} in $\mathbf{S^{n}}$ if it belongs to set of surfaces transverse to the each of these submanifolds. 

Thom's transversality theorem guarantees that the space of generic surfaces of degree $n$ is the union of everywhere dense open subsets of $\mathbf{S^{n}}$ (see \cite{ArnolSngI}, pg. 231).

If $S^{n}_{f}$ is generic, the \emph{set of hyperbolic points $\mathbb{H}_{f}$} consists of the generic hyperbolic, special hyperbolic, generic inflection and bi-inflection points. Moreover, the set $\mathbb{H}_{f}$ is an open subset of $S^{n}_{f}$. The \emph{set of elliptic points} of $S^{n}_{f}$ is denoted by $\mathbb{E}_{f}$, and is an open subset of $S^{n}_{f}$.

On the other hand, the \emph{set of parabolic points $\mathbb{P}_{f}$} of a surface generic $S^{n}_{f}$ consists of generic parabolic and special parabolic points. This set of points forms a smooth curve called \emph{parabolic curve}. This curve is the common topological boundary of the open sets of hyperbolic points and elliptic points. The special parabolic points are isolated on the parabolic curve.

On $\mathbb{H}_{f}$ there are two continuous direction fields, transverse each other. These direction fields are extended to the parabolic curve, where both coincide. We will referred to these direction fields  as \emph{asymptotic direction fields of $f$}. If  $S^{n}_{f}$ is generic, the unique asymptotic line is transverse at points of parabolic curve, except at special parabolic points, where the asymptotic line is tangent to the parabolic curve. In this case, the asymptotic direction fields of $f$ have well folded singularities at a special parabolic point: folded saddle, folded node and folded focus \cite{Dvdv}, \cite{Lk Dr}. 

The seven classes of points given above are invariant under affine and projective transformations. These classes were known at the end of the 19th century in the context of the enumerative geometry of complex algebraic surfaces, with works of Cayley, Zeuthen and Salmon, see \cite{Slmn} and \cite{Zthn}. The classification of points of a surface in $\mathbb{R}^{3}$ according to the order of tangency of their tangents is given by Landis E. E in \cite{Lnds}. The case of surfaces in $\mathbb{R}P^{3}$ was introduced by Platonova O. A. in \cite{Pltnv}. A classification in 22 types of surface points given by the action of the group of affine transformations of the three-dimensional space on the space of 3-jets of cooriented surfaces was introduced in \cite{Pnv1}.

\section{On the number of special parabolic points}

The special parabolic points are also known as \emph{Gaussian cusps} or \emph{cusp of the Gauss map} (see, for example \cite{Bnchff&Gffny&McCrry}) and \emph{godrons}. The last term is due to R, Thom \cite{Krgsn&thm}. For a generic surface each special parabolic point corresponds to a swallowtail point of the dual surface, that is, to an $A_{3}$ Legendre singularity (see \cite{Arnld4} and \cite{ArnolSngI}). In \cite{Urbcr}, Uribe-Vargas associates  to each special parabolic point a local projective invariant useful for the study of the local affine (projective) differential properties of swallowtails.  Bifurcations of the configuration formed by the parabolic and flecnodal curves that can occur in generic 1-parameter families of smooth surfaces is studied in \cite{Urbvltn}.

Finding the number of special parabolic points in generic algebraic surfaces, has been of special interest since the 19th century. In \cite{Klkv} is given the upper bound $2n(n-2)(11n-24)$ for the number of special parabolic points in generic algebraic surfaces of degree $n$ in $\mathbb{C}P^{3}$. This is also an upper bound of the number of special parabolic points on generic algebraic surfaces of degree $n$ in $\mathbb{R}P^{3}$. 

An upper bound, in the case of generic graphs of polynomial functions of degree $n$ in $\mathbb{R}^{3}$, is $(n-2)(5n-12)$, \cite{H&O&S1}. An upper bound for the number \(P_{+}\) (or \(P_{-}\)) of special parabolic points of index +1 (respectively -1) on a generic graph of polynomial function $f$ of degree $n$ in $\mathbb{R}^{3}$, such that its highest degree homogeneous polynomial has $k$ simple real linear factors (see \cite{Gdrrm&Ortz1}), are 
$$ P_{-} \leq \frac{(n-2)(8n-21) + k}{2} \text{ and } P_{+} \leq 1 + \frac{(n-2)(8n-21)-k}{2}.
$$ Until now, it is not known whether this bounds are optimal or not for $n\geq 4$.

If the set of non-hyperbolic points of a generic surface $S_{f}^{n}$ consists of at most $(2n-5)(n-3)+1$ mutually disjoint, compact, and convex disks, then all special parabolic points on $S_{f}^{n}$ have negative index (\cite{Urbcr}, and  \cite{H&O&S2}). Under these conditions an upper bound for the number of special parabolic points (having index $-1$) on $S_{f}^{n}$ is $$3(n-2)(n-3)+k,$$ where $k$ is the number of simple real linear factors  of the highest degree homogeneous polynomial of $f$ (see \cite{Gdrrm&Ortz1}). This upper bound is attained for $n=3,4$, but it is not known whether this bounds is optimal for $n\geq 5$.

\section{The Hessian curve of $f$}\label{hssncurve}

If a surface $S_{f}$ is the graph of a differentiable function $f$ on the xy-plane, the image of the parabolic curve of $S_{f}$ in the xy-plane under the projection $\pi(x, y, z)=(x, y)$ will be referred to as the \emph{Hessian curve of $f$}. It is defined by the zero level set of the function $$Hess(f)=f_{xx}f_{yy}-f_{xy}^{2}.$$ 

The directions determined by the projections of the asymptotic lines on the xy-plane under $\pi$ are the integral curves of the quadratic differential equation: 
$$II_{f} (x, y) = f_{xx}(x, y)dx^{2} + 2f_{xy}(x, y)dxdy + f_{yy}(x, y)dy^{2} = 0.$$ 

The binary form $II_{f} (x, y) $ will be referred to as the second fundamental form of $f$. For sake of simplicity we identify the solutions of this quadratic form with the
asymptotic directions of $f$. Also, we identify the set of hyperbolic points $\mathbb{H}_{f} \subset S_{f}$ with $\left\lbrace (x,y)\in \mathbb{R}^{2} \vert  Hess(f)(x,y)<0 \right\rbrace $, respectively $\mathbb{E}_{f} \subset S_{f}$ with the set $\left\lbrace (x,y)\in \mathbb{R}^{2} \vert  Hess(f)(x,y)>0 \right\rbrace $.

If $f\in \mathbb{R}[x,y]$ is a polynomial of degree $n\geq 3$, the Hessian curve of $f$ is a real affine curve of degree at most $2n-4$. When this curve is compact, there exists a unbounded connected component, denoted by $C_{u}$, contained in the complement of the Hessian curve in $\mathbb{R}^{2}$. In general, the Euler number of the projective closure $\hat{C_{u}}$ of the component $C_{u}$ is given by $\chi(\hat{C_{u}}) = 1- N_{u}$, where $N_{u}$ is the number of connected components of the boundary of $C_{u}$.

A basic aspect of the topology of a real algebraic curve is the  number of its connected components and the mutual positions in the (affine or projective) plane. This is connected to the first part of Hilbert's sixteenth problem. 

In 1876, Harnack proved in \cite{Hrnk} that a smooth plane curve in $\mathbb{R}P^{2}$ of even degree $d$ has at most $\frac{1}{2}(d-1)(d+2)+1$ connected components (usually called \emph{ovals}). This is also the upper bound for the number of connected components of a smooth compact curve  of  degree $d$ in $\mathbb{R}^{2}$. A non-compact curve in $\mathbb{R}^{2}$ of  degree $d$ can have at most $\frac{1}{2}(d-1)(d+2)$ ovals and $d$ unbounded components. These unbounded components come from the intersection of the corresponding projective curve with the line at infinity. 

By Harnack's Theorem, if the Hessian curve of a polynomial of degree $n$ is compact, then has at most $(2n-5)(n-3)+1$ ovals and, if it is a noncompact curve has, at most, $(2n-5)(n-3)$ ovals and $2n-4$ unbounded components. In \cite{Ortz}, Ortiz-Rodriguez exhibited a family of real polynomials of degree $n$ whose Hessian curve is smooth, compact and have $\frac{1}{2}(n-1)(n-2)$ ovals.  There are examples of compact Hessian curves of a polynomial $f$ of degree $n=3$ and $n=4$ with exactly $1$ and $4$ connected components. The first examples of compact Hessian curves of degree $4$ with four connected components were presented by   Ortiz-Rodr\'{i}guez and Sotille in \cite{Ort&Stll}. Using Viro patchworking in \cite{Brtrnd&Brgll}, Bertrand and Brugall\'{e}  construct a polynomial of degree $n\geq 4$ in two variables whose Hessian curve has $(n-4)^{2}$ ovals.  

The examples known for $n\geq5$ do not realize the Harnack bound. Indeed, in the compact case, the Harnack's bound for values $n = 5, 6, 7$ is $11, 22$ and $37$, respectively, whereas the number of connected components of Hessian curves of degree $n = 5, 6, 7$ that have been exhibited so far are $9, 16$ and $25$. Examples in small degrees of Hessian curve with unbounded components are exhibited in \cite{Ort&Stll}. 

The possible number and mutual position of the connected components of the Hessian curve in $\mathbb{R}^{2}$ of a polynomial $f\in \mathbb{R}[x,y]$ of degree $n\geq 3$ is concerned with problem 2000-1 in the list of Arnold's problems \cite{Arnld2b}, attributed to Ortiz-Rodr\'{i}guez. See also the related problems 2000-1, 2000-2, 2001-1, and 2002-1.
 
\section{The projective Hessian curve of $f$} \label{prhssncurve}

The projective closure of the Hessian curve of  $f\in \mathbb{R}[x,y]$ is given by $H_{f}(x,y,z)=0$, where $H_{f}$ is a homogeneous polynomial in $\mathbb{R}[x,y,z]$ such that $H_{f}(x,y,1) = Hess(f)$. This curve is called \emph{projective Hessian curve of} $f$. 

If we consider the homogeneous decomposition of $f=\sum _{i=r} ^{n} f_{i}$, where $f_{i}\in \mathbb{R}[x,y]$ is a homogeneous polynomial of degree $i$, then $$H_{f}(x,y,z) = \sum ^{2n-4}_{j=2r-4} h_{j}(x,y),$$ where $h_{2r-4} = Hess f_{2r-4}$ and $h_{2n-4} = Hess f_{n}$. 
In particular, the restriction of $H_{f}$ to the line at infinity $z = 0$ is $$H_{f}(x,y,0) = Hess f_{n}(x,y).$$
In case $f_{n}$ has no repeated real linear factors, then  $Hess f_{n}(x,y)$ is not identically zero (see \cite{Grvch}), and $H_{f}$ has degree $2n-4$.

In general, the projective Hessian curve of $f$ is a nonsingular curve of even degree equals $2n-4$. In this case, the complement in $\mathbb{R}P^{2}$ of the Hessian projective curve $H_{f}=0$ is the union of two disjoint open subsets, say $b^{-}$ and $b^{+}$. The set $b^{+}$ is an orientable smooth surface at which the sign of $H_{f}$ remains unchanged, while the open set $b^{-}$ is a nonorientable smooth surface at which $H_{f}$ takes the opposite sign. The topological closure in  $\mathbb{R}P^{2}$ of $b^{+}$ and $b^{-}$ will be denoted by $B^{+}$ and $B^{-}$, respectively.

An oval of a real algebraic curve in $\mathbb{R}P^{2}$ is called \emph{even} (\emph{odd}) if it is contained in an even (odd) number of ovals of the same curve. The number of even ovals is denoted by $P$ and the number of odd ovals by $N$.

These numbers $P$ and $N$ contain information about the topology of the surfaces $B^{+}$ and $B^{-}$. Indeed, the surface $B^{+}$ has $P$ connected components and the surface $B^{-}$ has $N+1$ connected components. According to Virgina Ragsdale in \cite{Rgsdl}, the Euler characteristics of these surfaces are $\chi(B^{+})=P-N$ and $\chi(B^{-})=N-P+1$. Furthermore, for any nonsingular real projective algebraic curve of even degree $m = 2k$, the number  $\chi(B^{+})=P-N$  satisfies the inequalities of Petrowsky \cite{Ptrwsk}: $$ -\frac{3}{2}k(k-1) \leq P-N \leq \frac{3}{2}k(k-1)+1.$$ 

In particular, if the Hessian projective curve $H_{f}=0$ has degree $2n-4$, we have $$ -\frac{3}{2}(n-2)(n-3) \leq P-N \leq \frac{3}{2}(n-2)(n-3)+1.$$ 

For our purposes, given a polynomial $f\in \mathbb{R}[x,y]$, it would be useful to determine the conditions under which the set $H_{\geq 0}$ coincides with surfaces $B^{+}$ or $B^{-}$. We will address this in more detail in the next section.

\section{Mutual position of the sets $\mathbb{H}_{f}$ and  $\mathbb{E}_{f}$}

We are interested in the basic affine geometric structure of a surface $S^{n}_{f}$. In particular, in this section we aim to describe some of topological features and the mutual position in the affine plane of the open sets  $\mathbb{H}_{f}$ and  $\mathbb{E}_{f}$, which  are imposed by the properties of $f_{n}$. The next homogeneous polynomial classes play a important role to fulfilled this goal.

A homogeneous polynomial on $\mathbb{R}[x,y]$ is called \emph{hyperbolic} (\emph{elliptic}) if its Hessian polynomial has no real linear factors and if it is non-positive (non-negative) at any point. 

An elliptic homogeneous polynomial in  $\mathbb{R}[x,y]$ has no real linear factors. On the other hand, a hyperbolic homogeneous polynomial has at least one real linear factor. Moreover, every real linear factor of a hyperbolic polynomial has multiplicity one \cite{Gdrrm&Ortz1}. In particular, every homogeneous polynomial $f$ of degree $n\geq 2$ with $n$ real linear factors, which are distinct up to non-zero constant factors, is  hyperbolic \cite{Gdrrm}.

Providing that $f_{n}$ is a hyperbolic or elliptic homogeneous polynomial, the Hessian curve of $f$ is compact and $C_{u}$ is a subset of  $\mathbb{H}_{f}$ or $\mathbb{E}_{f}$ (see \cite{Gdrrm&Ortz1}, \cite{H&O&S2}). Thus, if $f_{n}$ is hyperbolic, then $\mathbb{H}_{f}$ is a subset of the non-orientable surface $B^{-}$. Similarly, if $f_{n}$ elliptic, $\mathbb{E}_{f}$ is a subset of the non-orientable surface $B^{-}$. No other conditions on $f_{n}$ are known to determine whether the set of hyperbolic (or elliptic) points are orientable.

In \cite{ArnldAst}, Arnold stated the following conjecture regarding the connectedness of the space of homogeneous polynomial of degree $n$: 

{\it{ The number of connected components of the space of hyperbolic homogeneous polynomials of degree $n$ increases as $n$ increases (at least as a linear function of $n$)}}. 

In the same work, Arnold proves that the space of hyperbolic homogeneous polynomials of degree $3$ and $4$ are connected, but is disconnected for degree $n\geq 5$.

\section{Projection into the Poincar\'e sphere}
\label{sec:7}

For a polynomial $f\in\mathbb{R}\left[ x,y\right] $ of degree $n\geq 3$ let $\mathbb{X}_{1}$ and $\mathbb{X}_{2}$ be the two fields of asymptotic directions on the $xy$-plane defined by the equation  $II_{f}(x,y)=0$. Following ideas of Poincar\'{e} \cite{Pncr1}, 
we describe the extension of these direction fields to the line at infinity. Let $\mathbb{S}^{2} = \{(u,v,w)\in \mathbb{R}^3 \, | \, u^2+v^2+w^2=1\}$ be the unit sphere with center at the origin $\mathbf{0}$ in $\mathbb{R}^3$. Identify the tangent plane $T_{N}S^{2}$ of the sphere at the north pole $N=(0,0,1)$  with the plane xy-plane. For each point $\mathbf{x}=(x,y,1) \in T_{N}S^{2} $, the line through $\mathbf{0}$ and $\mathbf{x}$ intersects $\mathbb{S}^{2}$ at the points 
\begin{equation*}
s_{1}\left( \mathbf{x}\right) =\frac{1}{\sqrt{1+x^{2}+y^{2}}
} \left(
x,y,1\right), \,\,\,\, s_{2}\left( \mathbf{x}\right)
=-\frac{1}{\sqrt{1+x^{2}+y^{2}}} \left( x,y,1\right).
\end{equation*}

This defines two maps $s_{1}$ and $s_{2}$ from $\mathbb{R}^{2}$ to $\mathbb{S}^{2}$, called \emph{projections into the Poincar\'{e} sphere}. For $i=1,2$, the images, $s_{1}(\mathbb{X}_{i})$ and $s_{2}(\mathbb{X}_{i})$,  of the two fields of asymptotic directions $\mathbb{X}_{1}$ and $\mathbb{X}_{2}$ 
into both north and south hemispheres are the zero loci of the induced quadratic differential forms, $s_{1}^{\ast }\left( II_{f} \right)$
and $s_{2}^{\ast }\left( II_{f} \right)$, which are defined on the open north and south hemispheres, respectively. Moreover, the images of both fields into each open hemisphere consist of two fields of lines diffeomorphic to $\mathbb{X}_{1}$ and $\mathbb{X}_{2}$. 

As in \cite{Gnnz02}, the induced quadratic differential forms $s_{1}^{\ast }\left( II_{f}
\right)$ and $s_{2}^{\ast }\left( II_{f} \right)$ are extended to an analytical quadratic differential form  
\begin{equation}\label{EDLA}
\begin{pmatrix}
du & dv & d\omega%
\end{pmatrix}%
\begin{pmatrix}
\omega ^{2} F_{uu}\left( u,v,\omega \right) & \omega
^{2}F_{uv}\left( u,v,\omega \right) & \omega A\left(
u,v,\omega \right) \\
\omega ^{2}F_{uv}\left( u,v,\omega \right) & \omega
^{2} F_{vv}\left( u,v,\omega \right) & \omega B\left(
u,v,\omega \right) \\
\omega A\left( u,v,\omega \right) & \omega B\left( u,v,\omega \right) & S\left( u,v,\omega
\right)%
\end{pmatrix}%
\begin{pmatrix}
du \\
dv \\
d\omega%
\end{pmatrix}.%
\end{equation}
This quadratic forms is defined at each point in the sphere with the property that the equator is an integral curve of the fields defined by this form. Here, 
\begin{eqnarray*} 
&&F(u,v,\omega) = \sum_{i=0}^{n}\omega ^{n-i} f_i(u,v), \\ 
&&F_{uu}= \frac{\partial^2 F}{\partial u^2}, \,\, F_{uv}= \frac{\partial^2 F}
{\partial u\partial v}, \,\, F_{vv}= \frac{\partial^2 F}{\partial v^2},\\
&&A =-u F_{uu}-v F_{uv}, \\
&&B = -u F_{uv} -v F_{vv},\\
&&S = u^{2} F_{uu} +2uv F_{uv}+ v^{2} F_{vv}.
\end{eqnarray*}

The discriminant curve of (\ref{EDLA}) is the projective Hessian curve of $f$ together with the equator of the sphere $\mathbb{S}^{2}$. On the other hand, the quadratic differential form (\ref{EDLA}) has a \emph{singular point} at $(u,v,w)\in \mathbb{S}^{2}$ if {\rm (\ref{EDLA})} is zero at $(u,v,w)$. 

If $S_f$ is generic, all singular points of {\rm (\ref{EDLA})} (if there are any) are in the equator of the sphere $\mathbb{S}^{2}$. In fact, a point $\left( u_{0},v_{0},\omega _{0}\right)\in \mathbb{S}^{2}$, with
$\,\omega _{0}\neq 0$ is
a singular point of {\rm (\ref{EDLA})} if and only if the point $%
\left( x_{0},y_{0}\right) =\frac{1}{\omega _{0}}\left( u_{0},v_{0}\right) $ is a singular point of the second fundamental form $II_{f}(x,y)$.  A singular point $(u,v,w)\in \mathbb{S}^{2}$ of {\rm (\ref{EDLA})} is called \emph{singular point at infinity} if $w=0$. 

The number of singular points at infinity of (\ref{EDLA}) is always finite. This can be untrue for other binary quadratic differential equations (see \cite{Gnnz02}).

We denote by $\mathbb{Y}_{1}$ and $\mathbb{Y}_{2}$ the two fields of lines tangent to the sphere that are defined by the quadratic form (\ref{EDLA}). It is worth mentioning that these fields are not defined on the points $s_{1}\left( x,y,1 \right)$ and $s_{2}\left( x,y,1 \right)$ with $Hess(x,y)\geq 0$.  The set of singular points of {\rm (\ref{EDLA})} is also the set of singular points of the fields by $\mathbb{Y}_{1}$ and $\mathbb{Y}_{2}$. Moreover, if $f$ $\in\mathbb{R} \left[ x,y\right] $ is a polynomial of degree $n\geq 3$, then the set of singular points at infinity of the direction fields $\mathbb{Y}_i,\, i=1,2$  is  $$ \{(u,v,0) \in \mathbb{S}^{2} \, |\, f_{n}\left( u,v\right) =0\},
$$ where $\, f_n$ is the homogeneous part of degree $n$ of $f$. Therefore, the number of singular points at infinity of the direction fields $\mathbb{Y}_{i,}\, i=1,2$,  is twice the number of distinct real linear factors of the homogeneous polynomial $f_{n}$.

Let $p$ be a point on the equator of $\mathbb{S}^{2}$, such that $p$ is a singular point of the fields $\mathbb{Y}_{i,}\, i=1,2$. If $f_{n}$ has no repeated factors, then $H_{f} (p) < 0.$  Therefore, a singular point at infinity is locally surrounded by points in $H_{\leq 0}$.  In consequence, the Poincar\'e index of $\mathbb{Y}_{i,}\, i=1,2$, at $p$ is well defined.

\begin{figure}[htb]    
\begin{center}
\includegraphics[width=0.97in]{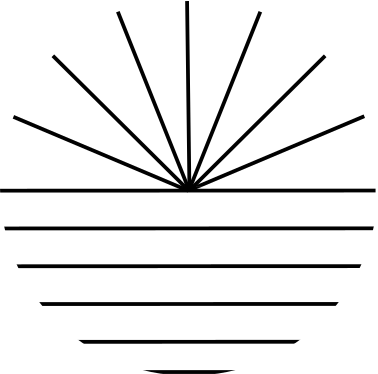}
\caption{Topological type of singular points at infinity.}
 \label{thetopotype}
\end{center}
\end{figure}

\begin{theorem}[\cite{Gdrrm&Ortz1}]\label{thindexonehalf}
Let $f = \sum_{i=1}^{n}f_{i} \in \mathbb{R}
\left[ x,y\right] $ be a polynomial of degree $n\geq 3.$ 
Assume that $f_{n}$ has no
repeated factors. Then, the Poincar\'{e} index of $\mathbb{Y}_{i}, i=1,2$ at each singular point at infinity is equal to $\frac{1}{2}$. Moreover, its topological type is shown in figure \ref{thetopotype}. 
\end{theorem}

For $i=1,2,\,$ Sing ($\mathbb{Y}_{i}$) denotes the set of singular points of the field $\mathbb{Y}_{i}$.

\begin{corollary}[\cite{Gdrrm&Ortz1}]\label{cotacasonoacotado}
Let $f$ $\in\mathbb{R}
\left[ x,y\right] $  be a polynomial of degree $n\geq 3$ such that $S_{f}^{n}$ is generic. Assume that $f_{n}$ has no repeated factors. If the projective Hessian curve of $f$  is transverse to the line at infinity, then 
\begin{equation*}
0\leq \sum_{\xi \in \mbox{Sing}\, (\mathbb{Y}_{i})} Ind_{\xi}\left( \mathbb{Y}_{i} \right) \leq n,
 \,\,\,\, \mbox{for} \,\, i=1,2. 
\end{equation*}
\end{corollary}

The restriction of the fields $\mathbb{Y}_{1}$ and $\mathbb{Y}_{2}$  
to the north hemisphere (or the south hemisphere) of $\mathbb{S}^{2}$ will be called {\it a projective extension of the two fields of asymptotic directions} $\mathbb{X}_1$ and $\mathbb{X}_2$. If $p$ is an singular point at infinity of the field $\mathbb{Y}_i$ there exists a 
projective extension on which $p$ is a singular point. We denote by $\tilde{\mathbb{X}}$ the projective extension that has the point $p$ as a singular point. Suppose that $f$ satisfies the assumptions of corollary 2 and that $f_{n}$ has $k$ linear real factors, then
\begin{equation} \label{sumptssing}
\sum_{\xi \in \mbox{Sing} (\tilde{\mathbb{X}})} Ind_{\xi}\left( \tilde{\mathbb{X}} \right) = \frac{k}{2}. 
\end{equation}
 
The next theorem provides a formula that relates the Euler 
characteristic of the surface $B^{\epsilon}$ with the Poincar\'e index of the singularities of fields $\tilde{\mathbb{X}}$ when it is defined on $B^{\epsilon}$, respectively. 

\begin{theorem}[\cite{Gdrrm&Ortz1}]\label{thglobalbminus}
Let $f$ $\in
\mathbb{R} \left[ x,y\right] $ be a polynomial of degree $n\geq 3,$ whose graph $S_f$ 
is generic and its projective Hessian curve is transverse to the line at infinity. 
Assume that $f_{n}$ has not repeated factors. If $\tilde{\mathbb{X}}$ is an extension to $\mathbb{R}P^2$ of a field of asymptotic directions. Then
\begin{equation}\label{maineqtn}
\sum_{\xi \in \mbox{Sing} (\tilde{\mathbb{X}})} Ind\left( \tilde{\mathbb{X}} \right) = \chi \left( B^{\epsilon}\right) +\frac{
	P_{-}-P_{+}}{2}\, 
\end{equation}
Where, $\epsilon$ is either $+$ or $-$, and $\tilde{\mathbb{X}}$ is defined on $B^{\epsilon}$. In both cases, $P_{-}$ denotes the number of special parabolic points with index $-1$ and $P_{+}$ the number of special parabolic points with index $+1$.
\end{theorem}

Other formulas involving the number of special parabolic points with index +1 and index -1 can be found in \cite{Blckr&Wlsn} and \cite{Bnchff&Thm}. In \cite{H&O&S2} is given an index formula for the field of asymptotic directions involving the number of connected components of the Hessian curve constituting the boundary of $C_{u}$, and the algebraic number of special parabolic points.

\section{Proof of the main Theorems}

\textbf{ Proof of theorem \ref{mainthrm1}}
 Let us suppose that there exists a polynomial  $f\in\mathbb{R} \left[ x,y\right] $ of degree $n$ with the properties stated in the Theorem \ref{mainthrm1}. 
Since $S_{f}$ is a generic surface such that the projective Hessian curve of $f$ is connected, transverse to the line at infinity and $B^{-}$ contains the set of hyperbolic points, then all the hypotheses of Theorem \ref{thglobalbminus} are satisfied. 

Therefore, the equality (\ref{maineqtn}) is satisfied.   Moreover $\chi(B^{-})=0$, since  $B^{-}$ is a closed M\"{o}bius band. Since $S_{f}$ has only one special parabolic point of positive index, then the number $P_{-}$ is equal to zero, and $P_{+}$ is equal to 1. From this we obtain that $$ \sum_{\xi \in \mbox{Sing} (\tilde{\mathbb{X}})} Ind_{\xi}\left( \tilde{\mathbb{X}} \right) = -\frac{1}{2}.$$ According to corollary \ref{cotacasonoacotado}, this equality is impossible since $$ \sum_{\xi \in \mbox{Sing} (\tilde{\mathbb{X}})} Ind_{\xi}\left( \tilde{\mathbb{X}} \right) \geq 0.$$

\textbf{ Proof of Corrollary \ref{maincor1}} 

Corrollary \ref{maincor1} follows from Theorem \ref{mainthrm1}, because  $C_{u}$ is contained in the set $H_{\leq 0}$.

\textbf{ Proof of theorem \ref{mainthrm2}}

Let us first consider the case where $B^{-}$ contains the set of hyperbolic points. Equation \ref{sumptssing} tell us that $ \sum_{\xi \in \mbox{Sing} (\tilde{\mathbb{X}})} Ind_{\xi}\left( \tilde{\mathbb{X}} \right) = \frac{k}{2}$, where $k$ is the number of simple real linear factors of $f_{n}$. Thus, from Theorem \ref{thglobalbminus}, we obtain that $$ \frac{k}{2}=  \chi \left( B^{-}\right) +\frac{P_{-}-P_{+}}{2}.$$  Therefore, if $0\leq k- 2\chi(B^{-})  $,  then $ P_{-}-P_{+}$ is greater than or equal to zero. 

The case where $B^{+}$ contains the set of hyperbolic points it follows in a similar way.

\end{document}